\documentclass{svmult}

\usepackage{makeidx}         
\usepackage{graphicx}        
\usepackage{multicol}        
\usepackage[bottom]{footmisc}
\usepackage{url}

\renewcommand{\PackageWarningNoLine}[2]{}
\usepackage{amsmath}
\usepackage{amsfonts}
\usepackage{color}
\pdfoutput=1

\begin{document}

\title*{A Direct Elliptic Solver Based on Hierarchically Low-rank Schur Complements}
\titlerunning{Hierarchically Low-rank Schur Complements}
\author{Gustavo Ch{\'a}vez, George Turkiyyah, \and David Keyes}
\authorrunning{G. Ch{\'a}vez, G. Turkiyyah \and D. Keyes}
\institute{\texttt{\{gustavo.chavezchavez,george.turkiyyah,david.keyes\}@kaust.edu.sa}, Extreme Computing Research Center, King Abdullah University of Science and Technology, Thuwal 23955, Saudi Arabia. (Submitted to the Proceedings of the 23th International Conference on Domain Decomposition Methods in Jeju, Korea)}

\maketitle

\abstract*{
A parallel fast direct solver for rank-compressible block tridiagonal linear systems is presented. 
Algorithmic synergies between Cyclic Reduction and Hierarchical matrix arithmetic operations result in a solver with $O(N \log^2 N)$ arithmetic complexity and $O(N \log N)$ memory footprint. 
We provide a baseline for performance and applicability by comparing with well known implementations of the $\mathcal{H}$-LU factorization and algebraic multigrid with a parallel implementation that leverages the concurrency features of the method.
Numerical experiments reveal that this method is comparable with other fast direct solvers based on Hierarchical Matrices such as $\mathcal{H}$-LU and that it can tackle problems where algebraic multigrid fails to converge.
}

\section{Introduction}

Cyclic reduction was conceived for the solution of tridiagonal linear systems,
such as the one-dimensional Poisson equation \cite{hockney65}. Generalized
to higher dimensions, it is known as block cyclic reduction (BCR) \cite{buzbee70}.
It can be used for general (block)
Toeplitz and (block) tridiagonal linear systems; however, it is 
not competitive for large problems, because its arithmetic complexity grows superlinearly. 
Cyclic reduction
can be thought of as a direct Gaussian elimination that recursively computes
the Schur complement of half of the system. Schur complement
computations have complexity that grows with the cost of the inverse, but by considering a tridiagonal system
and an even/odd ordering, cyclic reduction can decouple the system in such a way that the inverse of a
large block is the block-wise inverse of a collection of independent
smaller blocks. This addresses the
most expensive step of the Schur complement computation in terms of operation complexity
and does so in a way
that launches concurrent subproblems.
Its concurrency feature in the form of recursive bisection makes it interesting
for parallel environments, provided that
its arithmetic complexity can be improved. 

We address the time and memory complexity growth of the traditional cyclic reduction algorithm by approximating dense
blocks as they arise with hierarchical matrices ($\mathcal{H}$-matrices). The effectiveness of the block approximation 
relies on the rank structure of the original matrix.  
Many relevant operators are known to have low rank off the diagonal. This philosophy follows recent
work discussed below, but to our knowledge this is the first demonstration of the utility of complexity-reducing hierarchical substitution in the context of cyclic reduction.

The synergy of cyclic reduction and hierarchical matrices leads to a parallel fast direct solver of 
log-linear arithmetic complexity, $O(N \log^2 N)$, with controllable accuracy, using cyclic reduction outside and $\mathcal{H}$-Matrix arithmetic operations inside. The algorithm is 
purely algebraic, depending only on a block tridiagonal structure. Neither geometric information nor sophisticated re-orderings of the matrix are required. We call it 
Accelerated Cyclic Reduction (ACR).
Using a well-known implementation of $\mathcal{H}$-LU, we demonstrate the range of applicability of ACR over a set of model problems including strong convection and high frequency Helmholtz, including variable coefficient problems that cannot be tackled with the traditional FFT enabled version of cyclic reduction, FACR \cite{swarztrauber77}.  We show that it is competitive in time to solution with a global $\mathcal{H}$-LU factorization that does not exploit the cyclic reduction structure.
The fact that ACR is completely algebraic expands its range of applicability 
to problems with arbitrary coefficient structure within the block tridiagonal sparsity structure, subject to their amenability to rank compression. 
This gives the method robustness in some applications that are difficult for multigrid.
The concurrency and flexibility to tune the accuracy of individual matrix block approximations 
makes it interesting for emerging manycore architectures.  Finally, as with many direct solvers, there are complexity-accuracy tradeoffs that lead to the development of a new scalable preconditioner.

\section{Related Work}

Exploiting underlying low-rank structure is a trending strategy for improving the performance
of sparse direct solvers.

\begin{figure}
\centering
\includegraphics[height=3.1cm]{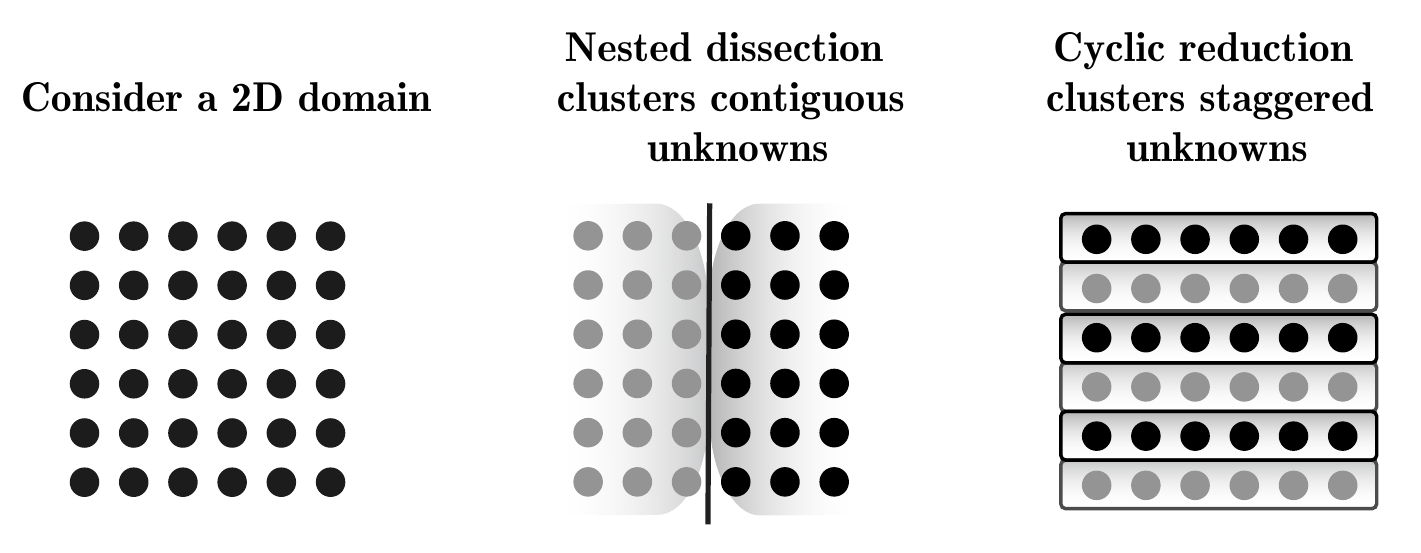}
\caption{The nested dissection ordering recursively clusters contiguous unknowns by bisection, whereas the
red/black ordering recursively clusters staggered unknowns, allowing isolation of a new readily manipulated diagonal block.}
\label{chavez_minisymposiumMS11_fig:acr_vs_nd}
\end{figure}

\textbf{Nested dissection based clustering of an $\mathcal{H}$-matrix} is known as
$\mathcal{H}$-Cholesky by Ibragimov {\em et al.} \cite{ibra07} and $\mathcal{H}$-LU by Grasedyck {\em et al.}
\cite{blackBoxHLU08,grasedyck09}, the main idea being to introduce $\mathcal{H}$-Matrix approximation on Schur complements based on domain decomposition. 
This is accomplished by a nested dissection ordering of the unknowns, and the advantage is that
large blocks of zeros are preserved after factorization. 
The non-zero blocks are replaced with low-rank approximations, and an LU factorization
is performed, substituting in $\mathcal{H}$ operations. 
Recently, Kriemann {\em et al.} \cite{Kriemann2014} demonstrated that $\mathcal{H}$-LU
implemented with a task-based scheduling based on a directed acyclic graph is well suited for 
modern many-core systems when compared with the conventional recursive algorithm. A similar line of work
by Xia {\em et al.} \cite{Xia10HSS_Cholesky} also proposes the construction of a rank-structured 
Cholesky factorization via the HSS hierarchical format. Figure \ref{chavez_minisymposiumMS11_fig:acr_vs_nd} illustrates the differences between nested dissection ordering and the even/odd (or red/black) ordering of cyclic reduction.

\textbf{Multifrontal factorization, with low-rank approximations of frontal matrices}, as in 
the work of Xia {\em et al.} \cite{xia10} also relies on nested dissection as the permutation strategy,
but it uses the multifrontal method as a solver. Frontal matrices are approximated with the HSS format,
while the solver relies on the corresponding HSS algorithms for elimination \cite{xia2010fast}. 
A similar line of work is the generalization of this method to 3D problems and general meshes by
Schmitz {\em et al.} \cite{Sch14,Sch12}. More recently, Ghysels {\em et al.} \cite{Ghysels15} introduced
a method based on a fast ULV decomposition \cite{chandrasekarana2006} and randomized
sampling of HSS matrices in a many-core environment, where HSS approximations are used to 
approximate fronts of large enough size, as the complexity constant in building an 
HSS approximation is only convenient for large matrices.

This strategy is not limited to any specific hierarchical format. Aminfar {\em et al.} \cite{AminfarD14}
proposed the use of the HODLR matrix format \cite{Ambikasaran2013},
also in the context of the multifrontal method. The well known solver MUMPS now also exploits 
the low-rank property of frontal matrices to accelerate its multifrontal implementation, as 
described in \cite{mumpsLR2014}.

\section{Accelerated Cyclic Reduction}

Consider the two-dimensional linear variable-coefficient Poisson equation (\ref{chavez_minisymposiumMS11_eq:probSta})
and its corresponding block tridiagonal matrix structure resulting from a second order 
finite difference discretization, as shown in (\ref{chavez_minisymposiumMS11_eq:system}):

\begin{equation}
\centering
-\nabla \cdot \kappa(\vec{x}) \nabla u=f(\vec{x}),
\label{chavez_minisymposiumMS11_eq:probSta}
\end{equation}

\begin{equation}
{\renewcommand{\arraystretch}{0.9}
A = \mbox{tridiag}(E_{i},D_{i},F_{i}) =
\begin{bmatrix}
~D_1 & F_1 &  &  &  \\
E_2 & D_2 & F_2 & & \\
 & ~\ddots & ~\ddots & ~\ddots & \\
 &  & E_{n-1} & D_{n-1} & F_{n-1} \\
 &  &  & E_{n} & D_{n}~\\
\end{bmatrix}.
}
\label{chavez_minisymposiumMS11_eq:system}
\end{equation}

We leverage the fact that for arbitrary $\kappa(\vec{x})$, the blocks $D_i$
are \emph{exactly} representable by rank 1 $\mathcal{H}$-Matrices, and the blocks $E_i$ and $F_i$ are diagonal.
As cyclic reduction progresses, the resulting blocks
will have a bounded increase in the numerical ranks of their off-diagonal blocks. 
This numerical off-diagonal rank may be tuned to accommodate 
for a specified accuracy. We choose the $\mathcal{H}$-Matrix
format proposed in \cite{hackbusch99} by Hackbusch, although ACR is not limited to a specific hierarchical
format. In terms of admissibility condition, we choose weak admissibility,
as the sparsity structure is known beforehand and it proved effective in our numerical experiments.

Approximating each block as an $\mathcal{H}$-matrix, we use
the corresponding hierarchical arithmetic operations as cyclic reduction progresses, instead of the conventional linear
algebra arithmetic operations. The following table summarizes the complexity estimates in 
terms of time and memory while dealing with a $n \times n$ block in a typical dense format
and as a block-wise approximation with a rank-$r$ $\mathcal{H}$-matrix.

\begin{center}
\begin{tabular}{c|c|c}
        & Inverse           & Storage\\ \hline
Dense Block     & $\mathcal{O}(n^3)$    & $\mathcal{O}(n^2)$ \\
$\mathcal{H}$ Block & $\mathcal{O}(r^2 n \log^2 n)$ & $\mathcal{O}(r n \log n)$
\label{chavez_minisymposiumMS11_table:advantages}
\end{tabular}
\end{center}

The following table summarizes the complexity estimates of the methods discussed so far in two dimensions, neglecting the dependence upon rank.

\begin{center}
  \begin{tabular}{c|c|c}
                  & Operations          & Memory \\ \hline
  BCR               & $\mathcal{O}(N^2)$          & $\mathcal{O}(N^{1.5} \log N)$ \\
  $\mathcal{H}$-LU  & $\mathcal{O}(N \log^2 N)$   & $\mathcal{O}(N \log N)$ \\
  ACR               & $\mathcal{O}(N \log^2 N)$   & $\mathcal{O}(N \log N)$ \\
  \end{tabular}
\label{chavez_minisymposiumMS11_table:summary}
\end{center}

With block-wise approximations in place, block cyclic reduction becomes ACR.
BCR consists of two phases: reduction and back-substitution. The reduction phase is equivalent
to block Gaussian elimination without pivoting on a permuted system $(PAP^T)(Pu)=Pf$. Permutation decouples the system,
and the computation of the Schur complement reduces the problem size by half.
This process is recursive and finishes when a single block is reached,
although the recursion can be stopped when the system is small enough to be solved directly. 
The second phase performs a conventional back-substitution to find the solution at every point of the domain.

As an illustration, consider a system of $n = 8$ points per dimension,
which translates into a $N \times N$ sparse matrix, with $N = n^2$. The first step
is to permute the system, which with an even/odd ordering becomes:

\begin{equation}
{\renewcommand{\arraystretch}{0.9}
\begin{array}{cccc}
\left[\begin{array}{cccc|cccc}
D_0 &     &     &     & F_0 &     &     &     \\
    & D_2 &     &     & E_2 & F_2 &     &     \\
    &     & D_4 &     &     & E_4 & F_4 &     \\
    &     &     & D_6 &     &     & E_6 & F_6 \\ \hline 
E_1 & F_1 &     &     & D_1 &     &     &     \\
    & E_3 & F_3 &     &     & D_3 &     &     \\
    &     & E_5 & F_5 &     &     & D_5 &     \\
    &     &     & E_7 &     &     &     & D_7 \\
\end{array}\right]
&
\left[\begin{array}{c}
u_0\\
u_2\\
u_4\\
u_6\\\hline 
u_1\\
u_3\\
u_5\\
u_7\\
\end{array}\right]
&
=
&
\left[\begin{array}{c}
f_0\\
f_2\\
f_4\\
f_6\\\hline 
f_1\\
f_3\\
f_5\\
f_7\\
\end{array}\right]
\end{array}.
}
\end{equation}

Consider the above $2\times2$ partitioned system as $H$.
The upper-left block is block-diagonal, which means that its inverse can be computed as the
inverse of each individual block ($D_0$, $D_2$, $D_4$, and $D_6$), 
in parallel and with hierarchical matrix arithmetics. The Schur complement of the upper-left partition may then be computed as follows:

\begin{equation}
{\renewcommand{\arraystretch}{0.9}
\begin{array}{cccc}
\left[\begin{array}{c|c}
H_{11}  & H_{12} \\ \hline 
H_{21}   & H_{22} \\\end{array}\right]
&
\left[\begin{array}{c}
u_{even} \\ \hline 
u_{odd} \\\end{array}\right]
&
=
&
\left[\begin{array}{c}
f_{even} \\ \hline 
f_{odd} \\\end{array}\right]
\end{array}.
\label{eq:sche}
}
\end{equation}

\begin{equation}
(H_{22}-H_{21}H_{11}^{-1}H_{12})u_{odd}=f^{(1)},\;\;\;\;\;\;\;\; f^{(1)} = f_{odd}-H_{21}H_{11}^{-1}f_{even}.
\label{eq:fund}
\end{equation}

Superscripts indicates algorithmic steps. A key property of the Schur complement
of a block tridiagonal matrix is that it also yields a block tridiagonal matrix, as can been seen in
the resulting permuted matrix system:

\begin{equation}
{\renewcommand{\arraystretch}{0.9}
\begin{array}{cccc}
\left[\begin{array}{cc|cc}
 D_0^{(1)}  &           & F_0^{(1)} &            \\
            & D_2^{(1)} & E_2^{(1)} & F_2^{(1)}  \\ \hline 
 E_1^{(1)}  & F_1^{(1)} & D_1^{(1)} &            \\
            & E_3^{(1)} &           & D_3^{(1)}  \\ 
\end{array}\right]
&
\left[\begin{array}{c}
u_0^{(1)}\\
u_2^{(1)}\\\hline 
u_1^{(1)}\\
u_3^{(1)}\\
\end{array}\right]
&
=
&
\left[\begin{array}{c}
f_0^{(1)}\\
f_2^{(1)}\\\hline 
f_1^{(1)}\\
f_3^{(1)}\\
\end{array}\right]
\end{array}.
}
\end{equation}

One step further, the computation of the Schur complement results in:

\begin{equation}
{\renewcommand{\arraystretch}{1.8}
\begin{array}{cccc}
\left[\begin{array}{cc}
D_0^{(2)} & F_0^{(2)}\\
E_1^{(2)} & D_1^{(2)}\\
\end{array}\right]
&
\left[\begin{array}{c}
u_0^{(2)}\\
u_1^{(2)}\\
\end{array}\right]
&
=
&
\left[\begin{array}{c}
f_0^{(2)}\\
f_1^{(2)}\\
\end{array}\right]
\end{array}.
}
\end{equation}

A last round of permutation and Schur complement computation
leads to the $D_0^{(3)}$ block, which is the last step of the
reduction phase of Cyclic Reduction. A conventional back-substitution
step recovers the solution.

\section{Numerical Results}

We select two test cases to provide a baseline of performance and robustness
as compared with high-performance implementations of $\mathcal{H}$-LU and AMG, namely HLIBpro \cite{HLib00},
and Hypre \cite{hypre-web-page}. Tests are performed in the shared memory environment of a 
12-core Intel Westemere processor.

\begin{figure}[h!]
\begin{minipage}[b]{0.47\linewidth}
\centering
\includegraphics[width=.95\linewidth]{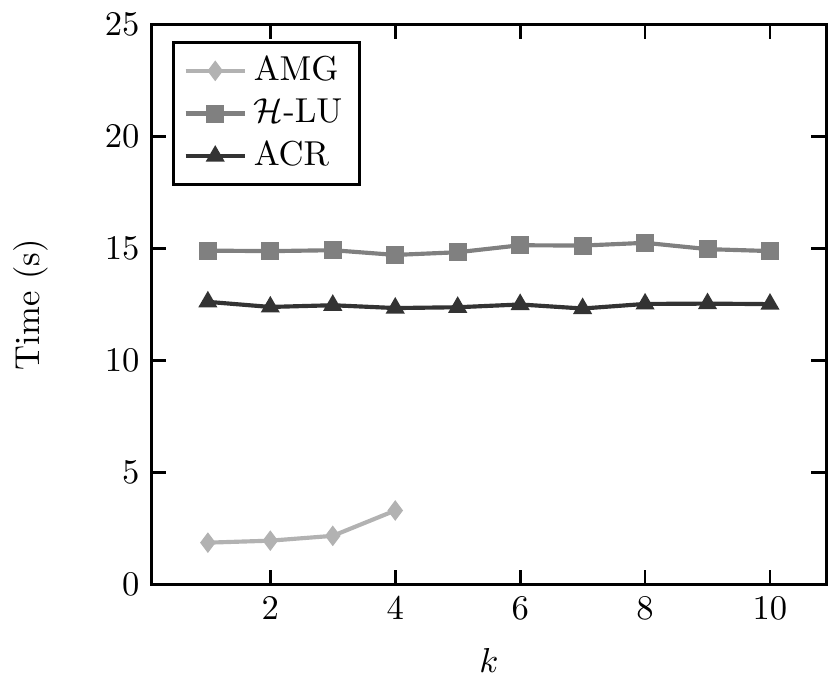}
\caption{Solvers performance while decreasing the number of points per wavelength. AMG fails to converge for large $k$.}
\label{chavez_minisymposiumMS11_figure3}
\end{minipage}
\hspace{0.5cm}
\begin{minipage}[b]{0.47\linewidth}
\centering
  \includegraphics[width=.95\linewidth]{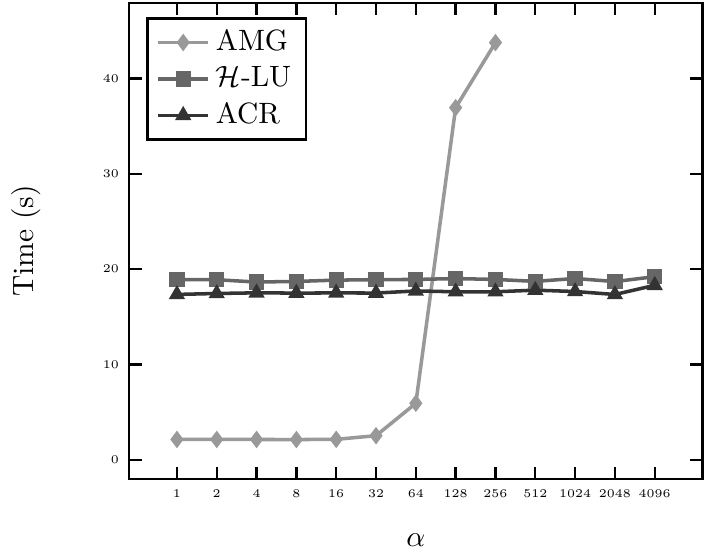}
  \caption{Stability of the solution as the convection dominance increases. AMG fails to converge for large $\alpha$.}
  \label{chavez_minisymposiumMS11_figure2}
\end{minipage}
\end{figure}

The first test is the wave Helmholtz equation. 
For large values of $kh$, where $h$ is the mesh spacing,
discretization leads to an indefinite matrix. Performance over a range of $k$ is
shown in
Figure \ref{chavez_minisymposiumMS11_figure3}, for $h=2^{-10}$.

\begin{equation}
\begin{aligned}
&\nabla^{2} u  + k^2u = f(\mathbf{x}), \;\;\;\;\;\;\;\;\;\;\;\; &\mathbf{x} \in \Omega = [0,1]^2  \;\;\; u(\mathbf{x}) = 0, \; x \in \Gamma &\\
&&f(\mathbf{x}) = 100 e^{-100((x-0.5)^2 + (y-0.5)^2)}.&\\
\end{aligned}
\end{equation}


The second test problem is the convection-diffusion equation with recirculating flow. The discretization of this equation leads to a nonsymmetric matrix.

\begin{equation}
\begin{aligned}
&-\nabla^{2} u  + \alpha b(\mathbf{x}) \cdot \nabla u = f(\mathbf{x}), \;\;\; &\mathbf{x} \in \Omega = [0,1]^2 \;\; u(\mathbf{x}) = 0, \; x \in \Gamma &\\
& b = {\cos(8\pi x) \choose \sin(8\pi y)} & f(\mathbf{x}) = 100 e^{-100((x-0.5)^2 + (y-0.5)^2)}. &
\end{aligned}
\end{equation}

\noindent
We progressively increase the convection dominance with $\alpha$,
accentuating the skew-symmetry. For small $\alpha$ AMG out performs ACR,
but since the increase in $\alpha$ does not depend on the rank structure of the matrix,
ACR maintains its performance for any $\alpha$. Performance can be seen in Figure 
\ref{chavez_minisymposiumMS11_figure2}

The discretization of 3D elliptic operators also
leads to a block tridiagonal structure, with the difference that each block is of size $n^2 \times n^2$,
as opposed to $n \times n$ as in the 2D discretization. Also, the numerical rank of the off-diagonal blocks 
grows faster than in the 2D case, which would lead into a superlinear solver of $\mathcal{O}(N^{\frac{4}{3}} \log N)$.
This complexity is not practical for 3D problems. Nonetheless, it is possible to use CR as a preconditioner,
as proposed in \cite{rodrigue84, reusken00}.
It is likewise advantageous to use ACR as a preconditioner since it is possible to
control the accuracy of the factorization, which directly translates into a tunable parameter of performance, 
while preserving the rich concurrency features of the method.
Preconditioning represents an interesting extension of the usability of ACR to tackle three dimensional
problems, with low-rank parameterization techniques. Further details are discussed in a companion paper in preparation with a larger set of test cases in two and three dimensions.

\section{Concluding Remarks}

We present a fast direct solver, ACR, for structured sparse linear systems 
that arise from the discretization of 2D elliptic operators. The solver approximates every block using an $\mathcal{H}$-Matrix, resulting in a log-linear
arithmetic complexity of $\mathcal{O}(N \log^2 N)$ with memory requirements of
$\mathcal{O}(N \log N)$. 

Robustness and applicability are demonstrated on model scalar problems and 
contrasted with established solvers based on the $\mathcal{H}$-LU factorization and algebraic multigrid.
Multigrid maintains superiority in scalar problems with sufficient definiteness and symmetry,
whereas hierarchical matrix based solvers can tackle problems where these properties are lacking.

Although being of the same asymptotic complexity as $\mathcal{H}$-LU, ACR has fundamentally different 
algorithmic roots which produce a novel alternative for a relevant class 
of problems with competitive performance, increasing concurrency as the problem grows,
and almost optimal memory requirements. In a companion paper we re-introduce the consideration of cyclic reduction as a preconditioner by exploiting tunable accuracy of low-rank approximation
and the exploitation of multicore features.

\bibliographystyle{plain}
\bibliography{chavez_minisymposiumMS11}

\end{document}